\documentclass[11pt,reqno,twoside]{article}
\usepackage{amsmath,amssymb,theorem,color,epsfig,epic,subfigure,amsfonts,hhline,graphicx}
\usepackage[normalem]{ulem}
\usepackage[colorinlistoftodos]{todonotes}
\theorembodyfont{\itshape}
\newtheorem{thm}{Theorem}[section]
\newtheorem{lem}[thm]{Lemma}
\newtheorem{prop}[thm]{Proposition}
{\theorembodyfont{\upshape}
\newtheorem{define}[thm]{Definition}
\newtheorem{rem}[thm]{Remark}
\newtheorem{ex}[thm]{Example}
}
\newtheorem{cor}[thm]{Corollary}

\newcommand{\Proof}[1][]{\noindent{\itshape Proof#1. }}
\newcommand{\EndProof}{\hfill$\Box$\\[1ex]}

\makeatletter
\let\@fnsymbol\@arabic
    
    \@addtoreset{equation}{section}
\makeatother

\def\hat{\widehat}

\def\sbs{\subset}

\def\a{\alpha}

   \def\L{\Lambda}             

\newcommand*{\Scale}[2][4]{\scalebox{#1}{$#2$}}%
\newcommand\unit{\hbox{\rm 1\kern-2.8truept l}}

\newcommand\Lform{{\mathcal{L}}\kern-7.56pt\raise1.0pt\hbox{$-$}}

\newcommand{\tr}{\mathrm{tr}}

\begin{document}
\title{ Conditional expectations associated with strongly quasi-invariant states and an application to spin systems}
\author{ Ameur Dhahri\footnote{Dipartimento di Matematica, Politecnico di Milano, Piazza Leonardo di Vinci 32, I-20133 Milano, Italy. E-mail: ameur.dhahri@polimi.it}, Chul Ki Ko\footnote{University College, Yonsei University, 85 Songdogwahak-ro,
Yeonsu-gu, Incheon 21983, Korea. E-mail: kochulki@yonsei.ac.kr}, and Hyun Jae Yoo\footnote{Department of  Applied Mathematics and Institute for Integrated Mathematical Sciences,
Hankyong National University, 327 Jungang-ro, Anseong-si,
Gyeonggi-do 17579, Korea. E-mail: yoohj@hknu.ac.kr}}
\date{ }
   \maketitle

\begin{abstract}
We discuss the conditional expectatins and martingales in relevance with $G$-strongly quasi-invariant states on a $C^*$-algebra $\mathcal A$, where $G$ is a separable locally compact group of $*$-automorphisms of $\mathcal A$. In the von Neumann algebra $\mathfrak A$ of the GNS representation, we define a unitary representation of the group and a group $\hat G$ of $*$-automorphisms of $\mathfrak A$, which is homomorphic to $G$. For the case of compact $G$,  we find a $\hat G$-invariant state on $\mathfrak A$ and define a conditional expectation with range the $\hat G$-fixed subalgebra. When $G $ is the union of increasing compact groups, we construct a sequence of conditional expectations and thereby construct (decreasing) martingales, which have limits by the martingale convergence theorem. As an example we consider $S_\infty$ the group of local permutations which acts on a $C^*$-algebra of infinite tensor product of finite dimensional $C^*$-algebras.  We also find an application in the classical spin systems.
\end{abstract}
\section{Introduction}
Since it was introduced by Umegaki \cite{Ume}, the conditional expectation for the operator algebra is an important concept not only just for an extension of classical concept but also for a tool to characterize the structure of the operator algebra. The concept has been further developed \cite{Acc-Ce, [Tak72]} and particularly it is a key concept to consider the martingales \cite{[Dang.Ngoc79], Hi, Tsu}.

Given a von Neumann algebra $\mathfrak A$ with a faithful normal state $\rho$, suppose that $(\mathfrak B_n)_n$ is an increasing (or decreasing, respectively) sequence of von Neumann subalgebras of $\mathfrak A$ and $(E_n)_n$ is a sequence of conditional expectations $E_n$ onto $\mathfrak B_n$ for each $n$. A sequence $(x_n)_n$ of $\mathfrak A$ is called martingale if (i) $x_n\in \mathfrak B_n$, (ii) $E_m(x_n)=x_m$ (respectively $E_n(x_m)=x_n$ in the decreasing case) for $m\le n$. The martingale convergence theorem says that if $(x_n)_n$ is a martingale, then there is an $x\in \mathfrak A$ such that $x_n\to x$ strongly as $n\to \infty$ \cite{Acc-Lo, [Dang.Ngoc79], Hi,Hi-Tsu,La,Tsu}.

In this paper we first investigate the martingale theory in relevance with $G$-strongly quasi-invariant states on a $C^*$-algebra $\mathcal A$. Here $G$ is a compact group or a separable locally compact group, which is the union of increasing compact groups of $*$-automorphisms of  $\mathcal A$. The main point lies in finding a state (weight) which is invariant under the group actions.

Let $\varphi$ be a faithful state on a $C^*$-algebra $\mathcal{A}$ and $G$ be a group of
$*$-automorphisms of $\mathcal{A}$. We say that $\varphi$ is $G$-strongly quasi-invariant on $\mathcal A$ if for all $g\in G$, there exists a self-adjoint operator $x_g\in\mathcal A$ such that  \cite{[Acc-Dha], Acc-Dha2}
$$\varphi(g(a))=\varphi(x_g a),\qquad \forall a\in \mathcal A.$$
In this case, it is proved in \cite{[Acc-Dha]} that for each $g\in G$, $x_g$ is strictly positive and it is an element of the centralizer of $\varphi$. Moreover, the algebra generated by the elements $x_g$ is a commutative $C^*$-algebra.

First, suppose $G$ is compact. In the von Neumann algebra of the representation of $(\mathcal A,\varphi)$, $\mathfrak A:=\pi(\mathcal A)''$, we define a unitary representation of $G$ and using this we induce a group $\hat G$, which is homomorphic to $G$, of $*$-automorphisms of $\mathfrak A$. Moreover, we find a $\hat G$-invariant state $\psi_G$ on $\mathfrak A$. Then we construct a conditional expectation $E:\mathfrak A\to \mathfrak B$, where $\mathfrak B$ is a von Neumann subalgebra of $\mathfrak A$ consisting of  $\hat G$-fixed elements.

When $G$ is a separable locally compact group, in particular $G$ has the structure $G=\bigcup_{N\in\mathbb{N}}G_{N}$,
where $(G_{N})_{N\in\mathbb{N}}$ is an increasing sequence of compact groups of
$*$-automorphisms of $\mathcal{A}$, imposing a condition that supports the existence of $\hat G$-invariant state, we construct a sequence of conditional expectations $(E_N)_N$ with ranges $\mathfrak B_N$'s, where for each $N$,  $\mathfrak B_N$ is a $\hat G_N$-invariant subalgebra of $\mathfrak A$. The sequence  $({E}_N)_N$ satisfies the martingale property, namely, for all $M\leq N$,
$${E}_N{E}_M(={E}_M{E}_N)={E}_N.$$
By this we have martingales: for any $x\in \mathfrak A$, the sequence $(x_n)_n$ with $x_n:=E_n(x)$ is a (decreasing) martingale. Then by a martingale convergence theorem such a sequence converges strongly to an element. This enables us to have a limit $E_N\to E_\infty$ as $N\to \infty$, where $E_\infty$ is an Umegaki conditional expectation onto $\mathrm{Fix}(u_G)$, the $\hat G$-invariant von Neumann subalgebra of $\mathfrak A$. Here the convergence is meant as a strong limit of $(E_N(x))_N$ to $E_\infty(x)$ for every $x\in \mathfrak A$ (Theorem \ref{thm:martingale_convergence}). We will apply these theory to nontrivial examples. One is the group of local permutations and the other one is the classical spin systems.

This paper is organized as follows. In Section \ref{sec:sqi_compact}, we briefly recall the definition of strongly quasi-invariant states with respect to compact groups. Next, given a $G$-strongly quasi-invariant state $\varphi$, we consider the GNS representation and induce a homomorphic group of automorphisms and an invariant state in the von Neumann algebra of the representation. Then, we construct a conditional expectation. In Section \ref{sec:Inductive_limit}, we discuss the inductive limit of compact groups. A sequence of conditional expectations will be considered. In Section \ref{sec:Martingale_convergence}, we discuss the martingales. Section \ref{sec:Example} is devoted to an example. We consider the  the locally compact group consisting of the local permutations on the set of nonnegative integers. This is the first example of a de Finetti theorem for quasi-invariant state in the noncommutative case.
In the final Section \ref{sec:spiin_system}, we apply the theory to the classical spin systems. We show that under certain conditions the Gibbs measures are $G$-strongly quasi-invariant, where $G$ is a locally compact group of spin flips or spin exchanges.
\section{Strongly quasi-invariant states with respect to a compact group}\label{subs-CG}\label{sec:sqi_compact}

In this section we consider the $G$-strongly quasi-invariant states for a compact group of $*$-automorphisms of a $C^*$-algebra and consider the conditional expectation.

\subsection{$G$-strongly quasi-invariant states}

Here we assume that $G$ is a compact group of  $*$-automorphisms of a $C^*$-algebra $\mathcal{A}$, and $\varphi$ is a $G$-strongly quasi-invariant faithful state on $\mathcal A$. We assume that the map $G\ni g\mapsto x_g\in \mathcal A$ is continuous. Denote by $\{\mathcal H,\pi,\Phi\}$ the cyclic representation of $(\mathcal A, \varphi)$. By \cite[Proposition 2.3.1]{BR}, the map $g\mapsto \pi(x_g)$ and consequently the map $g\mapsto \pi(x_g^{1/2})$ is continuous. From \cite{[Acc-Dha]}, it is proved that the map $U$ defined by
\begin{equation}\label{df--Ug}
U_g\pi(a)\Phi=\pi(g(a)x^{1/2}_{g^{-1}})\Phi\ ;\quad\forall\,a\in\mathcal{A},
\end{equation}
is a unitary representation of $G$ on $\mathcal H$. We assume the map $g\mapsto U_g$ is strongly continuous and on $\mathcal H$, define
$$P_G:=\int_GU_gdg,$$
where, and in the sequel, $dg$ denotes the normalized Haar measure on the compact group under consideration.
\begin{lem} \label{lem:range_projection}
$P_G$ is an orthogonal projection on $\mathcal H$ with range
\begin{eqnarray}\label{rang-p-g}
P_G(\mathcal H)=\{\xi\in \mathcal H:\;\;U_g(\xi)=\xi,\quad\forall g\in G\}=:\mathrm{Fix}_G(\mathcal H).
\end{eqnarray}
\end{lem}
\Proof We have
\begin{eqnarray*}
P_G^2&=&\int_G\left(\int_GU_gU_hdh\right)dg\\
&=&\int_G\Big(\int_GU_{gh}dh\Big)dg\\
&=&\int_G\Big(\int_GU_{h}dh\Big)dg\qquad(\mbox{The Haar measure is left translation invariant})\\
&=&\int_GP_Gdg=P_G.
\end{eqnarray*}
On the other hand, since the Haar measure is invariant by inversion, one gets
$$P_G^*=\int_GU_g^*dg=\int_GU_{g^{-1}}dg=\int_GU_{g}dg=P_G.$$
Now for any $\xi\in \mathcal H$ and $g\in G$,
\[
U_gP_G(\xi)=\left(\int_GU_gU_hdh\right)(\xi)=\left(\int_GU_{gh}dh\right)(\xi)=\left(\int_GU_{h}dh\right)(\xi)=P_G(\xi).
\]
So, $P_G(\mathcal H)\subset \text{Fix}_G(\mathcal H)$. Conversely, suppose $\xi\in \hbox{Fix}_G(\mathcal H)$. Then,
$$P_G(\xi)=\int_GU_g(\xi)dg=\int_G\xi dg=\xi.$$
Hence $\text{Fix}_G(\mathcal H)\subset P_G(\mathcal H)$. We conclude that $P_G$ is an orthogonal projection onto $\text{Fix}(G)$. 
\EndProof
From \eqref{df--Ug} it follows that
$$U_g\Phi=\pi(x_{g^{-1}}^{1/2})\Phi.$$
It holds that
\begin{equation}\label{def-K-G}
\Phi_G:=P_G\Phi=\Big(\int_GU_gdg\Big)\Phi=\int_G\pi\big(x_{g^{-1}}^{1/2}\big)dg\Phi=\int_G\pi\big(x_{g}^{1/2}\big)dg\Phi.
\end{equation}
Let us define the operator appearing in the r.h.s. of \eqref{def-K-G} by
\begin{equation}\label{eq:op_K}
K_G:=\int_G\pi\big(x_{g}^{1/2}\big)dg.
\end{equation}
The operator $K_G$ will play a central role in this paper and we emphasize here that the operators $P_G$ and $K_G$, both acting on $\mathcal H$, are not equal to each other in general but result in the same vector $\Phi_G$ when  applied to the vector $\Phi$. See Example \ref{ex:P_and_K} below.
Moreover, by \eqref{rang-p-g} $\Phi_G$ is $U(G)$-invariant where $U(G)=\{U_g:\, g\in G\}$:
\begin{equation}\label{eq:U(G)-invariance}
U_g\Phi_G=\Phi_G,\quad g\in G.
\end{equation}
From the continuity of the map $g\mapsto \pi(x_g^{1/2})$ one can show that the operator $K_G$ is bounded with a bounded inverse (see e.g. the proof of \cite[Theorem 1]{[Acc-Dha]}). Therefore,
\[
\Phi_G=P_G\Phi=K_G\Phi\neq 0.
\]
And for any $a\in \mathcal A$, one has
\begin{eqnarray}\label{n-r-phi}
\varphi(a)=\langle \Phi,\ \pi(a)\Phi\rangle=\langle K_G^{-1}\Phi_G,\ \pi(a)K_G^{-1}\Phi_G\rangle=\langle \Phi_G,\ K_G^{-1}\pi(a)K_G^{-1}\Phi_G\rangle.
\end{eqnarray}
Now define a state on the von Neumann algebra $\mathfrak A:=\pi(\mathcal A)''$ by
\begin{equation}\label{eq:G_inv_functional}
\psi_G(x):=\frac1{\|\Phi_G\|^2} \langle \Phi_G,x\Phi_G\rangle,\quad x\in \mathfrak A.
\end{equation}
For each $g\in G$, define a linear $*$-map on $\mathfrak A$ by
\begin{equation}\label{eq:inherited_group}
u_g(x):=U_g xU_g^*,\quad x\in \mathfrak A.
\end{equation}
In particular $u_g$ acts on $\pi(\mathcal A)$ as
\begin{equation}\label{eq:induced_group_action}
u_g(\pi(a))=\pi(g(a)), \quad a\in \mathcal A.
\end{equation}
In fact, by using the cocycle property of $x_g$'s we have $x_g^{-1}=g^{-1}(x_{g^{-1}})$ \cite{[Acc-Dha], Acc-Dha2}. Thus, for all $a,b\in \mathcal A$,
\begin{eqnarray*}
u_g(\pi(a))\pi(b)\Phi&=&U_g\pi(a)U_g^*\pi(b)\Phi\\
&=&U_g\pi(a)\pi(g^{-1}(b)x_g^{1/2})\Phi\\
&=&\pi\left(g(a)bg(x_g^{1/2})x_{g^{-1}}^{1/2}\right)\Phi\\
&=&\pi(g(a))\pi(b)\Phi.
\end{eqnarray*}
Therefore, by letting $\hat G:=\{u_g:g\in G\}$, we see that $\hat G$ is a group of $*$-automorphisms of $\mathfrak A$, which is homomorphic to $G$.
\begin{prop}\label{prop:G_invariance}
The state $\psi_G$ on $\mathfrak A$ defined in \eqref{eq:G_inv_functional} is $\hat G$-invariant.
\end{prop}
\Proof Take arbitrary $g\in G$ and $a\in \mathcal A$. Putting $\Psi_G(\cdot):=\|\Phi_G\|^2\psi_G(\cdot)$, by \eqref{eq:U(G)-invariance},
\begin{eqnarray*}
\Psi_G(u_g(\pi(a)))&=&\langle \Phi_G,U_g\pi(a)U_g^*\Phi_G\rangle\\
&=& \langle U_{g^{-1}}\Phi_G,\pi(a)U_{g^{-1}}\Phi_G\rangle\\
&=& \langle  \Phi_G,\pi(a) \Phi_G\rangle=\Psi_G(\pi(a)).
\end{eqnarray*}
Since $\pi(\mathcal A)$ is weakly dense in $\mathfrak A$, the proof is completed.
\EndProof
\begin{rem}\label{rem:G_invariant_state}
If we define a state $\varphi_G$ on $\mathcal A$ by
\begin{equation}\label{eq:G_invariant_state}
\varphi_G(a):=\psi_G(\pi(a)),\quad a\in \mathcal A,
\end{equation}
then by Proposition \ref{prop:G_invariance}, $\varphi_G$ is a $G$-invariant state. In fact,
\begin{eqnarray*}
\varphi_G(g(a))&=&\psi_G(\pi(g(a)))\\
&=&\psi_G(u_g(\pi(a)))\\
&=&\psi_G(\pi(a))\quad(\text{by Proposition \ref{prop:G_invariance}})\\
&=&\varphi_G(a).
\end{eqnarray*}
\end{rem}
\begin{ex}\label{ex:P_and_K}
Let us consider the Example 2 in \cite{DKY}. Let $\mathcal A=\mathcal M_2(\mathbb C)$.
Let $G=\{g_\theta|\theta=0, \frac{\pi}{2}, \pi, \frac{3\pi}{2}\}$, where
 \begin{equation}\label{eq:normal_state}
g_\theta(a):=R_{-\theta} aR_{\theta},\quad a\in \mathcal A,\quad R_\theta=\left(\begin{matrix}\cos\theta&-\sin\theta\\\sin\theta&\cos\theta\end{matrix}\right).
\end{equation}
Let
$$
\rho=\left(\begin{matrix}\frac{e^\beta}{1+e^\beta}&0\\0&\frac{1}{1+e^\beta}\end{matrix}\right).
$$
The state $\varphi(a):=\tr{(\rho a)}$, $a\in \mathcal A$, is $G$-strongly quasi-invariant with
 \begin{equation}\label{eq:RN_derivatives}
 x_{g_0}=x_{g_{\pi}}=I,\quad x_{g_{{\pi}/{2}}}=x_{g_{{3\pi}/{2}}}=\left(\begin{matrix}e^{-\beta}&0\\0&e^\beta\end{matrix}\right).
 \end{equation}
It is not hard to see that in the GNS representation $(\mathcal H,\pi,\Phi)$ of $(\mathcal A,\varphi)$ ($\mathcal H=\mathcal M_2(\mathbb C)$, $\Phi=I_2$, the $2\times 2$ unit matrix, and we use $\pi(a)=a$ for simplicity), the projection $P_G$ acts as
\begin{equation}\label{eq:projection}
P_G\pi(a)\Phi=\frac12\left(\begin{matrix}a_{11}+e^{-\beta/2}a_{22}&a_{12}-e^{\beta/2}a_{21}\\
a_{21}-e^{-\beta/2}a_{12}&a_{22}+e^{\beta/2}a_{11}\end{matrix}\right),\quad a=\left(\begin{matrix}a_{11}&a_{12}\\a_{21}&a_{22}\end{matrix}\right).
\end{equation}
Obviously, $P_G\pi(a)\Phi= \pi(Pa)$ with no $P\in \mathcal A$, meaning that $P_G\notin \pi(\mathcal A)$. On the other hand, one can directly show that
\begin{equation}\label{eq:op_K_ex}
K_G=\pi\left(\int_Gx_g^{1/2}dg\right)=\frac12\left(\begin{matrix}1+e^{-\beta/2}&0\\0&1+e^{\beta/2}\end{matrix}\right)\in \pi(\mathcal A).
\end{equation}
One sees, however, that the two operators $P_G$ and $K_G$ act on $\Phi=I_2$ resulting in
\[
P_G\Phi=K_G\Phi=\frac12\left(\begin{matrix}1+e^{-\beta/2}&0\\0&1+e^{\beta/2}\end{matrix}\right)\equiv \Phi_G\in \mathcal H.
\]
\end{ex}

\subsection{Conditional expectation}\label{subsec:conditional_expectation}
Let $G$ be a compact group of $*$-automorphisms of a $C^*$-algebra $\mathcal A$ and $\varphi$ a $G$-strongly quasi-invariant state on $\mathcal A$. We use the same notations introduced in the previous subsection. In particular, the von Neumann algebra $\mathfrak A=\pi(\mathcal A)''$ obtained by a GNS representation of $(\mathcal A,\varphi)$ is equipped with a $\hat G$-invariant state $\psi_G$. Define
\begin{equation}\label{eq:induced_group_invariant_subalgebra}
\text{Fix}(u_G):=\{x\in \mathfrak A:\,u_g(x)=x,\,\,\forall \,g\in G\}.
\end{equation}
Denoting $\mathfrak B:=\text{Fix}(u_G)$ a von Neumann subalgebra of $\mathfrak A$ we define a $*$-map $E:(\mathfrak A,\psi_G)\to \mathfrak B$ by
\begin{equation}\label{eq:conditional expectation}
E(x):=\int_Gu_g(x)dg,\quad x\in \mathfrak A.
\end{equation}
It is promptly shown that $E(x)\in \mathfrak B$. In fact, for any $g\in G$,
\[
u_g(E(x))=\int_G u_gu_h(x)dh=\int_G u_{gh}(x)dh=\int_G u_h(x)dh=E(x).
\]
To show that the map $E$ is an onto map, suppose that $x\in \mathfrak B$. Then, $u_g(x)=x$ for all $g\in G$ and hence $E(x)=\int_Gu_g(x)dg=\int_G xdg=x$ showing that $x$ is in the range of $E$.  
 \begin{thm}\label{thm:Umegaki_conditional_expectation}
 The map $E:(\mathfrak A,\psi_G)\to \mathfrak B$ is an Umegaki conditional expectation \cite{Ume}.
 \end{thm}
\Proof
 We have to show that $E$ is a normal contractive positive projection satisfying (i) $E(yxy')=yE(x)y'$ for all $x\in \mathfrak A$, $y,y'\in \mathfrak B$, (ii) $\psi_G\circ E=\psi_G$. Since the map $g\mapsto u_g$ is a $*$-automorphism and $G$ is compact, the normality, contractivity, and positivity are obvious. We check the remaining two properties. Let $x\in \mathfrak A$ and $y,y'\in \mathfrak B$. Then,
 \[
 E(yxy')=\int_Gu_g(yxy')dg=\int_Gu_g(y)u_g(x)u_g(y')dg=\int_Gyu_g(x)y'dg=yE(x)y'.
 \]
Also, for any $x\in \mathfrak A$, by the $\hat G$-invariance of $\psi_G$ (Proposition \ref{prop:G_invariance}),
\[
\psi_G(E(x))=\psi_G\left(\int_Gu_g(x)dg\right)=\int_G\psi_G(u_g(x))dg=\int_G\psi_G(x)dg=\psi_G(x).
\]
The proof is completed.  
\EndProof
\section{Inductive limit of compact groups}\label{sec:Inductive_limit}
Let $(G_{N})_{N\in\mathbb{N}}$ be an increasing sequence of compact groups of $*$-automorphisms
of $\mathcal{A}$ and let
\begin{equation}\label{separable}
G:=\bigcup_{N\in\mathbb{N}}G_{N},
\end{equation}
which is a locally compact group.  
Let $P_N\equiv P_{G_N}$
be the orthogonal projection onto $\hbox{Fix}_{G_N}(\mathcal H)$ as was defined in Section \ref{subs-CG}. Then for any $M\leq N$, $\hbox{Fix}_{G_N}(\mathcal H)\subset \hbox{Fix}_{G_M}(\mathcal H)$ and therefore $P_{N}\leq P_M$, in particular
\begin{equation}\label{marting}
P_NP_M=P_N.
\end{equation}
 Since $(P_N)_N$ is a decreasing sequence of orthogonal projections, it converges strongly to an orthogonal projection denoted by $P_G$:
$$P_G:=s-\lim_N P_N.$$
Denote by $\mathcal H_\infty$ the range of $P_G$ and recall that
$$\hbox{Fix}_G(\mathcal H):=\{\xi\in \mathcal H:\;\;U_g(\xi)=\xi,\quad\forall g\in G\}.$$
\begin{lem}\label{lemm2}  $\mathcal H_\infty=\mathrm{Fix}_{G}(\mathcal H)$.
\end{lem}
\Proof
 Notice by Lemma \ref{lem:range_projection}
\[
\mathcal H_\infty=\cap_N P_{N}(\mathcal H)=\cap_N \text{Fix}_{G_N}(\mathcal H).
\]
Thus the proof is completed by just noticing $\text{Fix}_G(\mathcal H)=\cap_N \text{Fix}_{G_N}(\mathcal H)$.
\EndProof
Denote by $\Phi_N:=P_N\Phi$, $\Phi_G:=P_G\Phi=\lim_NP_N\Phi=\lim_N K_N\Phi$, where $K_N:=K_{G_N}$ is given by \eqref{eq:op_K} for $G=G_N$. Throughout the paper we assume the Hypothesis (H) below.
\begin{enumerate}
\item[(H)] {\bf Hypothesis}. $\Phi_G\neq 0$.
\end{enumerate}
Below in Proposition \ref{prop:sufficient_condition} and Proposition \ref{prof:weak_limit}, we give sufficient conditions for (H). Under the hypothesis (H), we define a state $\psi_G$ on $\mathfrak A$ as in \eqref{eq:G_inv_functional}:  $\psi_G(\cdot)=\frac1{\|\Phi_G\|^2}\langle\Phi_G,\cdot\Phi_G\rangle $.
\begin{thm}\label{thm:G-invariant_functional} Suppose that the hypothesis (H) holds. Then the state $\psi_G$ is $\hat G$-invariant.
\end{thm}
\Proof
 Take any $a\in \mathcal A$ and $g\in G$. There exists an $N_g\in \mathbb N$
 such that $g\in G_N$ for all $N\ge N_g$. Then, for all $N\ge N_g$, by Proposition \ref{prop:G_invariance} we have
 \[
\langle \Phi_N,u_g(\pi(a))\Phi_N\rangle= \langle \Phi_N,\pi(a)\Phi_N\rangle.
\]
Taking the limit $N\to \infty$ in both sides, we get
\[
\langle \Phi_G,u_g(\pi(a))\Phi_G\rangle= \langle \Phi_G,\pi(a)\Phi_G\rangle,
\]
that is, $\psi_G(u_g(\pi(a)))=\psi_G(\pi(a))$, the $\hat G$-invariance of $\psi_G$.  
\EndProof
An immediate consequence of Theorem \ref{thm:G-invariant_functional} is the following.
\begin{cor}\label{cor:G-invariant_state} Under the Hypothesis (H), define a state $\varphi_G$ on $\mathcal A$ by
\[
\varphi_G(a):=\psi_G(\pi(a)),\quad a\in \mathcal A.
\]
Then, $\varphi_G$ is $G$-invariant.
\end{cor}
Let us now consider some sufficient conditions for the Hypothesis (H). Let $\lambda_N$ denote the normalized Haar measure on the group $G_N$ for each $N$.
\begin{prop}\label{prop:sufficient_condition}
Suppose that there are constants $\epsilon_0>0$, $\delta_0>0$, and $N_0\in \mathbb N$ such that for each $N\ge N_0$ there is a subset $A_N\subset G_N$ such that $\lambda_N(A_N)\ge \delta_0$ and for $g,h\in A_N$, $\varphi\big((x_gx_h)^{1/2}\big)\ge \epsilon_0$. Then, the Hypothesis (H) holds.
\end{prop}
\Proof
Recall that $x_g$'s commute with themselves. Let $N\ge N_0$. By the relation $\Phi_N=K_N\Phi$, we have
\begin{eqnarray}\label{eq:norm_bound_from_below}
\|\Phi_N\|^2&=&\langle K_N\Phi,K_N\Phi\rangle\nonumber\\
&=&\langle \Phi,\iint_{(G_N)^2}\pi\big((x_gx_h)^{1/2}\big)dhdg\Phi\rangle\nonumber\\
&=&\iint_{(G_N)^2}\varphi\big((x_gx_h)^{1/2}\big)dhdg\nonumber\\
&\ge&\iint_{(A_N)^2}\varphi\big((x_gx_h)^{1/2}\big)dhdg\nonumber\ge \epsilon_0(\delta_0)^2.
\end{eqnarray}
Here we have used the positivity of $x_gx_h$ since it is a product of commuting positive operators.  Since $\Phi_N\to \Phi_G$ strongly, we conclude $\|\Phi_G\|\ge \sqrt{\epsilon_0}\delta_0>0$.
\EndProof
Here we consider another sufficient condition for the Hypothesis (H).
\begin{prop}\label{prof:weak_limit}
Suppose that the sequence of operators $(K_N)_N$ converges weakly to an invertible operator $K_G$. Then, $\Phi_G$ is nonzero and $\psi_G$ is faithful. 
\end{prop}
\Proof
We have the relation
\begin{eqnarray}\label{eq:1st_step}
&&\|\Phi_N-K_G\Phi\|^2\nonumber\\
&=&\langle (K_N-K_G)\Phi,(K_N-K_G)\Phi\rangle\nonumber\\
&=&\langle K_N \Phi,K_N \Phi\rangle-\langle K_N \Phi, K_G\Phi\rangle-\langle K_G\Phi,K_N \Phi\rangle+\langle  K_G\Phi, K_G\Phi\rangle.
\end{eqnarray}
By the assumption, as $N\to \infty$, the second, third, and fourth terms together converges to $-\langle  K_G\Phi, K_G\Phi\rangle$. For the first term in the last line of \eqref{eq:1st_step} we see that
\begin{eqnarray*}
\lim_N\langle K_N \Phi,K_N \Phi\rangle&=&\lim_N \big(\langle \Phi_N-\Phi_G,\Phi_N\rangle+\langle \Phi_G,K_N\Phi\rangle\big)\\
&=&\langle \Phi_G,K_G\Phi\rangle\\ 
&=&\lim_N \langle K_N\Phi,K_G\Phi\rangle\\
&=&\langle K_G\Phi,K_G\Phi\rangle
\end{eqnarray*}
We have shown that $\Phi_G=K_G\Phi$ and it is nonzero since $K_G$ is invertible. Particularly, it also implies that  $\Phi_G$ is a cyclic vector for $\mathfrak A$. By \cite[Proposition 4.1]{DKY}, $\Phi_g:=\pi(\sqrt{x_g})\Phi$ belongs to the positive cone $\mathcal P$ associated with $\Phi$. Obviously, $\Phi_N=\int_{G_N}\Phi_gdg$ also belongs to $\mathcal P$ and therefore $\Phi_G$ is an element of $\mathcal P$. Now, by \cite[Proposition 2.5.30]{BR}, $\Phi_G$ is also separating for $\mathfrak A$.  Therefore, $\psi_G$ is faithful. 
\EndProof
\section{Martingales}\label{sec:Martingale_convergence}

We continue with a setting of the previous section. $G=\cup_NG_N$ is a separable locally compact group consisting of an increasing sequence of compact groups of $*$-automorphisms of a $C^*$-algebra $\mathcal A$. $\varphi$ is a $G$-strongly quasi-invariant state on $\mathcal A$. We have a von Neumann algebra $\mathfrak A=\pi(\mathcal A)''$ of the GNS representation of $(\mathcal A,\varphi)$. We assume the Hypothesis (H). Therefore, $\mathfrak A$ is equipped with a $\hat G$-invariant state $\psi_G(\cdot)=\frac1{\|\Phi_G\|^2}\langle \Phi_G,\cdot\Phi_G\rangle$. For each $N$, define
\begin{equation}\label{eq:Umegaki}
{E}_N(x):=\int_{G_N}u_g(x)dg,\quad x\in \mathfrak A.
\end{equation}
By Theorem \ref{thm:Umegaki_conditional_expectation}, $E_N$ is an Umegaki conditional expectation onto $\mathfrak B_N:=\text{Fix}(u_{G_N})$. In particular, the following relations hold:
\begin{eqnarray}\label{def-widetilde-EN}
&{E}_N(\pi(a))=\int_{G_N}u_g(\pi(a))dg =\int_{G_N}\pi(g(a))dg,\quad\forall a\in \mathcal A,&\\
&u_h{E}_N(\pi(a))=\int_{G_N}\pi(hg(a))dg={E}_N(\pi(a)) \label{invariance}&
\end{eqnarray}
Obviously, $(\mathfrak B_N)_N$ is a decreasing sequence of von Neumann subalgebras of $\mathfrak A$.
\begin{thm}\label{thm:martingale}
Suppose the Hypothesis (H). It holds that $\psi_G\circ E_N=\psi_G$ for all $N$ and the sequence $(E_N)_N$ satisfies the martingale property, namely, for $M\le N$,
\begin{equation}\label{eq:martingale}
E_NE_M(=E_ME_N)=E_N.
\end{equation}
\end{thm}
\Proof
By Theorem \ref{thm:G-invariant_functional}, the state $\psi_G=\frac1{\|\Phi_G\|^2}\langle\Phi_G,\cdot \Phi_G\rangle$ is $u_g$-invariant for all $g\in G$. Therefore, for all $x\in \mathfrak  A$ and $N\in \mathbb N$,
\begin{eqnarray*}
\psi_G( E_N(x))&=&\psi_G\left(\int_{G_N}u_g(x)dg\right)\\
&=&\int_{G_N}\psi_G(u_g(x))dg\\
&=&\int_{G_N}\psi_G(x)dg\\
&=&\psi_G(x).
\end{eqnarray*}
Now suppose that $M\le N$. One has
\begin{eqnarray*}
{E}_N{E}_M(\pi(a)) &=&E_N\left(\int_{G_M}\pi(h(a))dh\right)\\ 
 &=&\int_{G_M}\left(\int_{G_N} \pi(gh(a))dg\right)dh\\
  &=&\int_{G_M}\left(\int_{G_N} \pi(g(a))dg\right)dh\\
&=&\int_{G_M}{E}_N(\pi(a))\,dh={E}_N(\pi(a)).
\end{eqnarray*}
The relation $E_ME_N=E_N$ can be similarly shown and the proof is completed.
\EndProof
\begin{thm}\label{thm:martingale_convergence}
Suppose that the Hypothesis (H) is satisfied and the state $\psi_G$ is faithful, or suppose that the sequence $(K_N)_N$ converges weakly to an invertible operator $K_G$. Then, for any $x\in \mathfrak A$, the sequence $(x_N)_N$, where $x_N:=E_N(x)$, is a (decreasing) martingale, and hence has a strong limit. Furthermore, by defining $E_\infty(x):=\lim_NE_N(x)$, $E_\infty$ is an Umegaki conditional expectation onto $\mathrm{Fix}(u_G)$.
\end{thm}
\Proof
Since we are in the decreasing situation, we have to show for $M\le N$, $E_N(x_M)=x_N$ \cite{[Dang.Ngoc79], Hi, Hi-Tsu,La,Tsu}. For any $y_N\in \mathfrak B_N$ we compute
\begin{eqnarray*}
\psi_G(E_N(x_M)y_N)&=&\psi_G(E_N(x_My_N))\\
&=&\psi_G(x_My_N)\\
&=&\psi_G(E_M(x)y_N)\\
&=&\psi_G(E_M(xy_N))\\
&=&\psi_G(xy_N)\\
&=&\psi_G(E_N(xy_N))\\
&=&\psi_G(E_N(x)y_N)\\
&=&\psi_G(x_Ny_N).
\end{eqnarray*}
We have shown that $\psi_G((E_N(x_M)-x_N)y_N)=0$ for all $y_N\in \mathfrak B_N$. Taking $y_N=(E_N(x_M)-x_N)^*$ and using the faithfulness we conclude that $E_N(x_M)-x_N=0$. The martingale convergence theorem is well known \cite{[Dang.Ngoc79], Hi, Hi-Tsu,La,Tsu}. Now let us define 
\begin{equation}\label{eq:martingale_limit}
E_\infty(x):=\lim_NE_N(x),\quad x\in \mathfrak A. 
\end{equation}
Obviously, $\text{Ran}E_\infty\subset \cap_N\text{Fix}(u_G)=\text{Fix}(u_G)$. On the other hand, suppose that $x\in \text{Fix}(u_G)$. Then, for each $N\in \mathbb N$, $x\in \text{Fix}(u_{G_N})$ and $E_N(x)=x$. Thus $E_\infty(x)=\lim_NE_N(x)=x$ showing that $E_\infty:\mathfrak A\to\text{Fix}(u_G)$ is an onto map. Furthermore, for any $y,y'\in \text{Fix}(u_G)$ and $x\in \mathfrak A$, we have 
\[
E_\infty(yxy')=\lim_NE_N(yxy')\overset{y,y'\in \text{Fix}(u_{G_N})}=\lim_N(yE_N(x)y')=yE_\infty(x) y'.
\]
And, for any $x\in \mathfrak A$,
\[
\psi_G\circ E_\infty(x)=\psi_G(\lim_NE_N(x))=\lim_N\psi_G\circ E_N(x)\overset{\psi_G\circ E_N=\psi_G}=\psi_G(x).
\]
Therefore, $E_\infty:(\mathfrak A,\psi_G)\to \mathfrak B_\infty\equiv \mathrm{Fix}(u_G)$ is an Umegaki conditional expectation.
\EndProof
 
\section{Example: the group of local permutations}\label{sec:Example}
Let $\mathbb N$ be the set of nonnegative integers and let $\mathcal S_\infty=\cup_{N\in\mathbb N}\mathcal S_N$ be the group of local permutations on $\mathbb N$. Putting $\mathcal B:=\mathcal B(\mathbb C^d)$, let $\mathcal{A}=\otimes_{n\in \mathbb{N}}\mathcal B$ be the $C^*$-algebra of infinite direct product of copies $\mathcal B$ \cite{Ta}. To say more precisely, for each finite set $F\subset \mathbb N$, let $\mathcal A_F:=\otimes_{n\in F}\mathcal B$ be the $C^*$-algebra of finite direct product of $\mathcal B$'s. For $F_1\subset F_2$, there is a natural embedding $\mathcal A_{F_1}\subset \mathcal A_{F_2}$ and $\mathcal A=\otimes_{n\in \mathbb{N}}\mathcal B$ is defined as the $C^*$-inductive limit of $\mathcal A_F$'s. We may consider $\mathcal A_F$ as a $C^*$-subalgebra of $\mathcal A$ and in particular for each $n\in \mathbb N$, $\mathcal A_n\equiv j_n(B):=(\otimes_{k\le n-1}\unit_\mathcal B)\otimes \mathcal B\otimes(\otimes_{k\ge n+1}\unit_B)$ is a subalgebra obtained by embedding $\mathcal B$ into the $n$th position.

For each $n\in \mathbb N$, let $W_n\in \mathcal B$ be a density matrix and let $\varphi(\cdot) :=\otimes_{n\in\mathbb{N}} \tr\left(W_n\cdot\right)$ be a
product state on $\mathcal{A}$. We assume  $[W_n,W_m]=0$ for all $m,n\in\mathbb N$. In a moment we will see that $\varphi$ is $G$-strongly quasi invariant for $G=\mathcal S_\infty$, but here we remark that we are considering the finite dimensional algebra  $\mathcal B(\mathbb C^d)$, otherwise we need to consider the generalized $G$-strongly quasi invariance \cite{DKY}.
Let $(b_n)_n$ be a sequence of element of $\mathcal B$ such that $b_n=\unit$ except finitely many  $n$'s. Then
$$\varphi(\prod_{n\in \mathbb N}j_n(b_n))=\prod_{n\in \mathbb N}\varphi_n(b_n)\;''=''\;\tr_{\mathbb N}\Big(\prod_{n\in \mathbb N}j_n(W_n)\prod_{n\in\mathbb N}j_n(b_n)\Big)$$
where $\varphi_n(b_n)=\tr(W_n b_n)$ is a state on $\mathcal B$ and $\tr_{\mathbb N}(\cdot)=\otimes_{\mathbb N}\tr(\cdot)$. Then for any $\sigma\in\mathcal S_N$, one has
 \begin{eqnarray*}
\varphi\Big(\sigma\Big( \prod_nj_n(b_n)\Big) \Big)&=&\tr_{\mathbb N}\Big(\prod_nj_n(W_n)\sigma\Big(\prod_nj_n(b_n) \Big)\Big)\\
&=&\tr_{\mathbb N}\Big(\sigma^{-1}\Big(\prod_nj_n(W_n)\Big)\prod_nj_n(b_n)\Big)\\
&=&\tr_{\mathbb N}\Big(\prod_nj_n(W_{\sigma(n)})\Big)\prod_nj_n(b_n)\Big)\\
&=&\tr_{\mathbb N}\Big(\prod_nj_n(W_n)\prod_nj_n(W_{\sigma(n)}/W_n) \prod_nj_n(b_n)\Big).
\end{eqnarray*}
Therefore, one gets
\begin{equation}\label{eq:RN_derivative}
x_\sigma=\prod_{n\in\Lambda_\sigma}j_n(W_{\sigma(n)}/W_n),
\end{equation}
where $\Lambda_\sigma$ is the support of $\sigma$ meaning that $\sigma(j)=j$ for $j\notin\Lambda_\sigma$.

\begin{figure}[h]
\begin{center}
\includegraphics[width=0.90\textwidth]{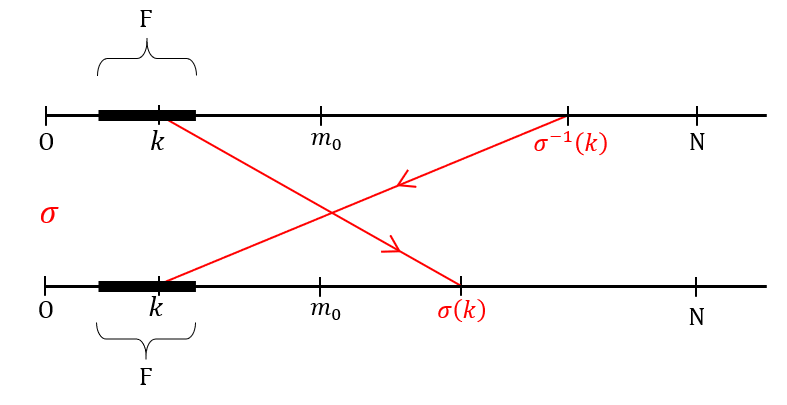}
\caption{ $\sigma\in \mathcal S_N^{(m_0)}$}\label{fig:outer_permutation}
\end{center}
\end{figure}

Let $W\in \mathcal B$ be a fixed density matrix and let $F$ be a finite subset of $\mathbb N$. We assume that
\begin{equation}\label{eq:finitely_generated}
W_n=W,\quad\forall n\in F^c.
\end{equation}
We also assume that there exists a constant $C>1$ such that
\begin{equation}\label{eq:uniform_bound}
\frac1C\le \|W_n\|\le C,\quad n\in \mathbb N.
\end{equation}
The Hypothesis (H) was crucial in this paper. We first check that (H) holds in this model.
\begin{lem}\label{lem:(H)}
The Hypothesis (H) holds for the above model.
\end{lem}
\Proof
For the proof we will use Proposition \ref{prop:sufficient_condition}. Fix an $m_0\in \mathbb N$ such that $F\subset \{0,1,\cdots,m_0\}$. For $N>m_0$, define (see Figure \ref{fig:outer_permutation})
\begin{equation}\label{eq:outer}
\mathcal S_N^{(m_0)}=\{\sigma\in \mathcal S_N:\,\sigma(k)>m_0,\,\,\sigma^{-1}(k)>m_0\text{ if }k\in \{0,1,\cdots,m_0\}\},
\end{equation}
and we put $A_N:=\mathcal S_N^{(m_0)}$. Denoting by $|A|$ the cardinality of a set $A$, it is not hard to see that
\begin{equation}\label{eq:asymptotic_diminishing}
\lim_{N\to \infty}\frac{1}{N!}|\mathcal S_N\setminus\mathcal S_N^{(m_0)}|=0.
\end{equation}
Therefore, we can find a $\delta_0>0$ and $N_0>m_0$ such that for all $N\ge N_0$, $\lambda_N(A_N)\ge \delta_0$.
Now, for $N\ge N_0$, if $\sigma\in A_N$, by \eqref{eq:RN_derivative} we see that
\begin{eqnarray*}
x_\sigma&=&\prod_{n\in\Lambda_\sigma}j_n(W_{\sigma(n)}/W_n)\\
&=&\prod_{n\in F}j_n(W/W_n)j_{\sigma^{-1}(n)}(W_n/W).
\end{eqnarray*}
Therefore, for $\sigma,\tau\in A_N$ we compute
\begin{eqnarray*}
\varphi((x_\sigma x_\tau)^{1/2})&=&\prod_{\substack{n\in F:\\\sigma^{-1}(n)\neq \tau^{-1}(n)}}\tr(W_{\sigma^{-1}(n)}^{1/2}W^{1/2})\tr(W_{\tau^{-1}(n)}^{1/2}W^{1/2})\\
&\ge& C^{-2 N(\sigma,\tau)}\ge C^{-2|F|}>0,
\end{eqnarray*}
where $N(\sigma,\tau)=|\{n\in F: \sigma^{-1}(n)\neq \tau^{-1}(n)\}|$.
Thus, the conditions of Proposition \ref{prop:sufficient_condition} are fulfilled and we are done.
\EndProof
In the next proposition we show a stronger result, namely, the weak convergence of $(K_N)$. It is another support for the Hypothesis (H) by Proposition \ref{prof:weak_limit}.
\begin{prop}\label{prop:example}
Under the assumptions \eqref{eq:finitely_generated} and \eqref{eq:uniform_bound}, the sequence $(K_N)_N$ converges weakly to an operator $K_G$ given by
\begin{equation}\label{eq:K_G}
K_G=\left(\prod_{n\in F}\tr(W_n^{1/2}W^{1/2})\right)\pi\left(\prod_{n\in F}j_n\left((W/W_n)^{1/2}\right)\right).
\end{equation}
\end{prop}
\Proof
Fix a number $n_0\in\mathbb N$ such that $F\subset\{0,1,\cdots,n_0\}$. Let $a=\prod_{n\in \mathbb N}j_n(a_n)$ and $b=\prod_{n\in \mathbb N}j_n(b_n)$ be the elements of $\mathcal A$ such that there exists a $k_0\in \mathbb N$ and $a_n=\unit$ and $b_n=\unit$ for $n\ge k_0$. For big enough $N$'s such that $N>m_0:=\max\{n_0,k_0\}$ we decompose $\mathcal S_N$ as
\[
\mathcal S_N=\mathcal S_N^{(m_0)}\cup(\mathcal S_N\setminus\mathcal S_N^{(m_0)}),
\]
where $\mathcal S_N^{(m_0)}$ is defined in \eqref{eq:outer}.
In order to show the weak convergence of $(K_N)_N$, let us compute the limit of $\langle \pi(a)\Phi,K_N\pi(b)\Phi\rangle$. From \eqref{eq:RN_derivative}, we have (we will omit the representation symbol ``$\pi$'' whenever there is no danger of confusion)
\[
K_N=\frac1{N!}\sum_{\sigma\in \mathcal S_N}\prod_{n\in \Lambda_\sigma}j_n(W_{\sigma(n)}^{1/2}/W_n^{1/2}).
\]
Therefore,
\begin{eqnarray*}
&&\langle \pi(a)\Phi,K_N\pi(b)\Phi\rangle\\
&=&\frac1{N!}\sum_{\sigma\in \mathcal S_N}\tr_{\mathbb N}\left(\Big(\prod_{n\in \mathbb N}j_n(W_n)\Big)\Big(\prod_{n\in \mathbb N}j_n(a_n^*)\Big)\Big(\prod_{n\in \Lambda_\sigma}j_n(W_{\sigma(n)}^{1/2}/W_n^{1/2})\Big)\Big(\prod_{n\in \mathbb N}j_n(b_n)\Big)\right)\\
&=&\Scale[0.80]{\frac1{N!}\left(\sum_{\sigma\in \mathcal S_N^{(m_0)}}+\sum_{\sigma\in \mathcal S_N\setminus\mathcal S_N^{(m_0)}}\right)\tr_{\mathbb N}\left(\Big(\prod_{n\in \mathbb N}j_n(W_n)\Big)\Big(\prod_{n\in \mathbb N}j_n(a_n^*)\Big)\Big(\prod_{n\in \Lambda_\sigma}j_n(W_{\sigma(n)}^{1/2}/W_n^{1/2})\Big)\Big(\prod_{n\in \mathbb N}j_n(b_n)\Big)\right)}\\
&=:& (\text{outer})_N+(\text{inner})_N.
\end{eqnarray*}
By using \eqref{eq:asymptotic_diminishing}, it can be easily shown that
\begin{equation}\label{eq:inner}
|(\text{inner})_N|\le \frac{1}{N!}|\mathcal S_N\setminus \mathcal S_N^{(m_0)}|\cdot C^{4|F|}M^{2k_0}\underset{N\to \infty}\longrightarrow 0,
\end{equation}
where $M:=\max\{\|a_n\|,\,\|b_n\|:\,n\in \mathbb N\}$. In order to estimate the first term, suppose that $\sigma\in \mathcal S_N^{(m_0)}$. Then, one sees that (see Figure \ref{fig:outer_permutation})
\begin{eqnarray*}
\prod_{n\in \Lambda_\sigma}j_n(W_{\sigma(n)}^{1/2}/W_n^{1/2})&=&\left(\prod_{k\in F}j_k(W_{\sigma(k)}^{1/2}/W_k^{1/2})\right)\left(\prod_{k\in F}j_{\sigma^{-1}(k)}(W_{k}^{1/2}/W_{\sigma^{-1}(k)}^{1/2})\right)\\
&=&\left(\prod_{k\in F}j_k(W^{1/2}/W_k^{1/2})\right)\left(\prod_{k\in F}j_{\sigma^{-1}(k)}(W_k^{1/2}/W^{1/2})\right).
\end{eqnarray*}
Also for those $\sigma$, we see that $a_k=\unit$ and $b_k=\unit$ for $k\in \sigma^{-1}(F)$.
Therefore, if $\sigma\in \mathcal S_N^{(m_0)}$, then
\begin{eqnarray*}
&&\tr_{\mathbb N}\left(\Big(\prod_{n\in \mathbb N}j_n(W_n)\Big)\Big(\prod_{n\in \mathbb N}j_n(a_n^*)\Big)\Big(\prod_{n\in \Lambda_\sigma}j_n(W_{\sigma(n)}^{1/2}/W_n^{1/2})\Big)\Big(\prod_{n\in \mathbb N}j_n(b_n)\Big)\right)\\
&=&\tr_{\mathbb N}\left(\Big(\prod_{n\in \mathbb N}j_n(W_n)\Big)\Big(\prod_{n\in \mathbb N}j_n(a_n^*)\Big)K_G\Big(\prod_{n\in \mathbb N}j_n(b_n)\Big)\right)\\
&=&\langle\pi(a)\Phi,K_G\pi(b)\Phi\rangle,
\end{eqnarray*}
where $K_G$ is defined in \eqref{eq:K_G}.
Now since $\frac1{N!}|\mathcal S_N^{(m_0)}|\underset{N\to \infty}\longrightarrow 1 $, we have
\begin{equation}\label{eq:outer}
\lim_{N\to \infty} (\text{outer})_N=\langle\pi(a)\Phi,K_G\pi(b)\Phi\rangle.
\end{equation}
Combining \eqref{eq:inner} and \eqref{eq:outer} we see that $K_n$ converges weakly to $K_G$.
\EndProof
\begin{rem}\label{rem:ex_group_invariance}
Let us compute the $G$-invariant state $\varphi_G$ in Corollary \ref{cor:G-invariant_state} for the example in this section. Let $a=\otimes_{n\in \mathbb N}a_n$ be any element of the form in the beginning of the proof of Proposition \ref{prop:example}. By noticing $\|\Phi_G\|^2=\|K_G\Phi\|^2=\left(\prod_{n\in F}\tr(W_n^{1/2}W^{1/2})\right)^2$,
\begin{eqnarray*}
&&\varphi_G(a)\\
&=&\psi_G(\pi(a))\\
&=&\frac1{\|\Phi_G\|^2}\langle \Phi_G,\pi(a)\Phi_G\rangle\\
&=&\frac1{\|K_G\Phi\|^2}\langle K_G\Phi,\pi(a)K_G\Phi\rangle\\
&=&   \tr_{\mathbb N}\left(\Big(\prod_{n\in \mathbb N}j_n(W_n)\Big)\Big(\prod_{n\in F}j_n((W/W_n)^{1/2})\Big)\Big(\prod_{n\in \mathbb N}j_n(a_n)\Big)\Big(\prod_{n\in F}j_n((W/W_n)^{1/2})\Big)\right)\\
&=&  \tr_{\mathbb N}\left(\Big(\prod_{n\in F}j_n(W_n)\Big)\Big(\prod_{n\in F^c}j_n(W)\Big)\Big(\prod_{n\in F}j_n((W/W_n)^{1/2})\Big)^2\Big(\prod_{n\in \mathbb N}j_n(a_n)\Big) \right)\\
&=& \tr_{\mathbb N}\left( \Big(\prod_{n\in \mathbb N}j_n(W)\Big)\Big(\prod_{n\in \mathbb N}j_n(a_n)\Big) \right)\\
&=& \prod_{n\in\mathbb N}\tr(Wa_n).
\end{eqnarray*}
Obviously, $\varphi_G$ is $G$-invariant and is the infinite product of the identical states: $\varphi_G=\otimes_{n\in \mathbb N}\varphi_n$, $\varphi_n(\cdot)=\tr(W\cdot)$.
\end{rem}
\section{An application to classical spin systems}\label{sec:spiin_system}
In this section we consider the classical spin systems in the statistical mechanical models. 
\subsection{Gibbs measures}
Let $\mathbb Z^d$ be the $d$-dimensional integer space. Let $\Omega_0=\{1,-1\}$ be the set representing  the spins at each site. $\Omega_0$ is equipped with a Bernoulli distribution $\mu_0$: $\mu_0(\{1\})=\mu_0(\{-1\})=\frac12$. We let $\mu_\Lambda$ be the probability measure $\mu_0^\Lambda$ on the set $\Omega_\Lambda:=\Omega_0^\Lambda$.  The whole configuration space is denoted by  $\Omega:=\Omega_{\mathbb Z^d}$. For each $\Lambda\sbs \mathbb Z^d$, $\mathcal F_\Lambda$ denotes a $\sigma$-algebra on the set $\Omega$ generated by the spin variables in the set $\Lambda$. We simply denote by $\mathcal F$ for $\mathcal F_{\mathbb Z^d}$. In the sequel, when $\Lambda\sbs \mathbb Z^d$ is a finite subset we  denote it by $\Lambda\sbs\sbs\mathbb Z^d$.

We abuse the notations but following the tradition of statistical mechanics, by an interaction $\Phi=(\Phi_\Lambda)_{\Lambda\sbs\sbs\mathbb Z^d}$ we mean a set of real-valued functions $\Phi_\Lambda:\Omega\to \mathbb R$, which is $\mathcal F_\Lambda$-measurable. 
\begin{ex}[Ising model]\label{ex:Ising}
The interaction for Ising model with a nearest neighborhood interaction is given by for $\xi=(\xi_i)_{i\in \mathbb Z^d}\in \Omega$
\[
\Phi_\Lambda(\xi)=\begin{cases}-J\xi_i\xi_j,&\Lambda=\{i,j\},\,\,|i-j|=1,\\ h\xi_i,&\Lambda=\{i\},\\
0,&\text{otherwise}\end{cases}.
\]
Here $J$ is the interaction strength and $h$ denotes the external magnetic field strength; $J>0$ for ferro magnetic model and $J<0$ for anti-ferro magnetic model. 
\end{ex}
In the sequel we assume that the interaction is translation invariant, i.e., $\Phi_X(\sigma_i(\omega))=\Phi_{\sigma_i^{-1}(X)}(\omega)$, where $\sigma_i$ is the translation by $i$, i.e., for all $j\in \mathbb Z^d$, $\sigma_i(j)=j+i$ and $(\sigma_i(\omega))_j=\omega_{j+i}$ for $\omega\in \Omega$. Furthermore, we assume that 
\begin{equation}\label{eq:interaction_condition}
\sum_{X\ni 0}\|\Phi_X\|<\infty.
\end{equation}
Given an interaction $\Phi$, for each $\Lambda\sbs\sbs \mathbb Z^d$  define a function $H_\Lambda(\cdot|\cdot)$, the Hamiltonian with boundary condition, by
\begin{equation}\label{eq:Hamiltonian}
H_\Lambda^\Phi(\zeta|\omega):=\sum_{X\cap\Lambda\neq\emptyset}\Phi_\Lambda(\zeta_\Lambda\omega_{\Lambda^c}),\quad \zeta,\,\,\omega\in \Omega,
\end{equation}
where  $\zeta_\Lambda\omega_{\Lambda^c}\in \Omega$ is the juxtaposition of $ \zeta_\Lambda$ and $\omega_{\Lambda^c}$, which are the restrictions of $\zeta$ and $\omega$ on $\Lambda$ and $\Lambda^c$, respectively. Let us define 
\[
Z_\Lambda^\Phi(\omega):=\int_{\Omega_\Lambda}\exp\left[-H_\Lambda^\Phi(\zeta|\omega)\right]d\mu_\Lambda(\zeta_\Lambda).
\]
We give a definition of a Gibbs measure for the interaction $\Phi$. \cite{Ge}.
\begin{define}\label{def:Gibbsian_specification}
Let $\Phi$ be an interaction satisfying \eqref{eq:interaction_condition} and $\omega\in \Omega$, $\Lambda\sbs\sbs\mathbb Z^d$. Then the probability measure 
\[
A\mapsto \gamma_\Lambda^\Phi(A|\omega):=\frac1{Z_\Lambda^\Phi(\omega)}\int_{\Omega_\Lambda}\exp\left[-H_\Lambda^\Phi(\zeta|\omega)\right]1_{A}(\zeta_\Lambda\omega_{\Lambda^c})d\mu_\Lambda(\zeta_\Lambda)
\]
on $(\Omega,\mathcal F)$ is called the Gibbs distribution in $\Lambda$ with boundary condition $\omega$ and interaction  $\Phi$. The system $(\gamma_\Lambda^\Phi)_{\Lambda\sbs\sbs\mathbb Z^d}$ is called the Gibbsian specification for $\Phi$. Any probability measure $\mu$ on $(\Omega,\mathcal F)$ is called a Gibbs measure for the interaction $\Phi$ if it satisfies the so called Dobrushin-Lanford-Ruelle (DLR) equations:
\begin{equation}\label{eq:DLR}
\mathbb E_\mu[A|\mathcal F_{\Lambda^c}]=\gamma_\Lambda^\Phi(A|\cdot) \,\,\mu\text{-a.s. }\text{for all }A\in \mathcal F\text{ and }\Lambda\sbs\sbs\mathbb Z^d.
\end{equation}
The set of Gibbs measures for the interaction $\Phi$ is denoted by $\mathcal G(\Phi)$.
\end{define}
\begin{rem}\label{rem:existence}
(i) It is known that under the condition \eqref{eq:interaction_condition} $\mathcal G(\Phi)\neq \emptyset$ \cite{Ge,Pr}.\\
(ii) For any continuous function $f\in C(\Omega)$ on $\Omega$, define
\[
\gamma_\Lambda^\Phi(f|\omega):=\int_\Omega f(\zeta_\Lambda\omega_{\Lambda^c})\gamma_\Lambda^\Phi(d\zeta|\omega)=\frac1{Z_\Lambda^\Phi(\omega)}\int_{\Omega_\Lambda}\exp\left[-H_\Lambda^\Phi(\zeta|\omega)\right]f(\zeta_\Lambda\omega_{\Lambda^c})d\mu_\Lambda(\zeta_\Lambda).
\]
The DLR equation \eqref{eq:DLR} is equivalent to saying that 
\begin{equation}\label{eq:DLR2}
\mu(f):=\int f(\omega)d\mu(\omega)=\mu(\gamma_\Lambda^\Phi(f|\cdot)),\quad f\in C(\Omega),\,\,\L\sbs\sbs\mathbb Z^d.
\end{equation}
\end{rem}
\subsection{Group actions and strong quasi-invariance of Gibbs measures}
From now on we assume that the interaction is of finite range, i.e.,  there is an $R>0$ such that $\Phi_X=0$ if $\text{diam}(X)>R$.

Notice that $\Omega$ is a compact space and let $\mathcal A$ be the $C^*$-algebra $C(\Omega)$ equipped with the sup-norm. Any probability measure on $(\Omega,\mathcal F)$, in particular any Gibbs measure, is a state on $\mathcal A$. For each $N\in \mathbb N$, let $\Lambda_N=[-N,N]^d\cap \mathbb Z^d$ denote the rectangular box with sides of length $2N+1$. Let $(G_N)_{N\in \mathbb N}$ be an increasing sequence of automorphisms of $\mathcal A$ such that for each $N\in \mathbb N$, $G_N$ depends only on the local configurations in $\Lambda_N$. We let $G=\cup_{N\in \mathbb N}G_N$. For the group $G$, mostly we have in mind the group of spin interchanges or spin flips defined as follows:
\begin{ex}\label{ex:spin_automorphism_group}
We consider the following group of automorphisms. Notice that any continuous bijection $\tau:\Omega\to \Omega$ naturally induces an automorphism $\tau:\mathcal A\to \mathcal A$ by 
\[
\tau(f)(\omega)=f(\tau(\omega)),\quad f\in \mathcal A.
\]
\begin{enumerate} 
\item[(i)]  (Spin exchanges) For $i\neq j\in \mathbb Z^d$, $\tau_{ij}:\Omega\to \Omega$ is defined by
\[(\tau_{ij}(\omega))_k=\omega_k^{ij}:=\begin{cases}\omega_k,&k\neq i,j\\ \omega_j,&k=i,\\
\omega_i, &k=j.\end{cases}. 
\]
\end{enumerate}
The group $G_N$ is generated by $\{\tau_{ij}:i\neq j\in \Lambda_N\}$. In other words, $G_N$ consists of spin permutations in the box $\Lambda_N$.
\item[(ii)] (Spin flips) For each $i\in \mathbb Z^d$, $\tau_i:\Omega\to \Omega$ is defined by
\[(\tau_i(\omega))_j=\omega_j^i:=\begin{cases}\omega_j,&j\neq i,\\ -\omega_i,&j=i.\end{cases}.  
\]
The group $G_N$ is generated by $\{\tau_i:i\in \Lambda_N\}$ so it is the group of partial spin flips in $\Lambda_N$. 
\end{ex}
\begin{thm}\label{thm:quasi_invariance_spin_classical}
Let $\Phi$ be an interaction satisfying \eqref{eq:interaction_condition} and let $\mu$ be a Gibbs measure for $\Phi$. Let $G=\cup_{N\in \mathbb N}G_N$ be one of the locally compact groups introduced in Example \ref{ex:spin_automorphism_group}. Then $\mu$ is $G$-strongly quasi-invariant with cocycles $x_\tau$ given by 
\begin{equation}\label{eq:cocycles_spin}
x_\tau(\omega)=\exp[H(\omega)-H(\tau^{-1}(\omega))], \quad \tau\in G,
\end{equation}
here, and in the sequel, the exponent is defined by 
\[
H(\omega)-H(\tau^{-1}(\omega))=\lim_{N\to \infty} \sum_{X\sbs \Lambda_N}\left(\Phi_X(\omega)-\Phi_X(\tau^{-1}(\omega))\right)\Big], 
\]
which is well-defined since $\tau$ gives only a local change. 
\end{thm} 
 \Proof
Suppose that $\tau\in G_{N_0}$. Let $N>N_0+R$ and fix $\Lambda\sbs\sbs\mathbb Z^d$ such that $\Lambda_N\sbs \Lambda$. By \eqref{eq:DLR2}, we have
\begin{eqnarray*}
\mu(\tau(f))&=&\mu(\gamma_{\Lambda_N}^\Phi(\tau(f)|\cdot))\\
&=&\int\left(\frac1{Z_\Lambda^\Phi(\omega)}\int_{\Omega_\Lambda}\exp\left[-H_\Lambda^\Phi(\zeta|\omega)\right]\tau(f)(\zeta_\Lambda\omega_{\Lambda^c})d\mu_\Lambda(\zeta_\Lambda) \right)d\mu(\omega) \\
&=&\int\left(\frac1{Z_\Lambda^\Phi(\omega)}\int_{\Omega_\Lambda}\exp\left[-H_\Lambda^\Phi(\zeta|\omega)\right]f(\tau(\zeta_\Lambda)\omega_{\Lambda^c})d\mu_\Lambda(\zeta_\Lambda) \right)d\mu(\omega) \\
&=&\int\left(\frac1{Z_\Lambda^\Phi(\omega)}\int_{\Omega_\Lambda}\exp\left[-H_\Lambda^\Phi(\tau^{-1}(\zeta)|\omega)\right]f(\zeta_\Lambda\omega_{\Lambda^c})d\mu_\Lambda(\zeta_\Lambda) \right)d\mu(\omega) \\
&=&\int\left(\frac1{Z_\Lambda^\Phi(\omega)}\int_{\Omega_\Lambda}\exp\left[-H_\Lambda^\Phi(\zeta|\omega)\right](x_\tau f)(\zeta_\Lambda\omega_{\Lambda^c})d\mu_\Lambda(\zeta_\Lambda) \right)d\mu(\omega),
\end{eqnarray*}
where 
\begin{eqnarray*}
x_\tau(\zeta)&=&\exp\Big[\sum_{X\cap \Lambda_{N_0}\neq \emptyset}\left(\Phi_X(\zeta)-\Phi_X(\tau^{-1}(\zeta))\right)\Big]\\
&=&\exp[H(\zeta)-H(\tau^{-1}(\zeta))].
\end{eqnarray*}
\EndProof
Let $(\mathcal H_\mu,\pi_\mu,\Phi_\mu)$ be the GNS representation of  $(\mathcal A,\mu)$. We can think of $\mathcal H_\mu=L^2(\Omega,\mu)$, $\pi_\mu(f)=f$ for $f\in C(\Omega)$, $\mathfrak A=\pi_\mu(\mathcal A)''=L^\infty(\Omega,\mu)$ acting as multiplication operators on $L^2(\Omega,\mu)$, and $\Phi_\mu=1$, the unit function on $\Omega$. Let us compute the unitary operators $U_\tau$, $\tau\in G_N$, in \eqref{df--Ug}. 

Given an $f\in C(\Omega)$ and $\tau\in G$, recall from \eqref{df--Ug} that
\begin{equation}\label{eq:U_again}
U_\tau\pi_\mu(f)\Phi_\mu=\pi_\mu(\tau(f)x_{\tau^{-1}}^{1/2})\Phi_\mu,
\end{equation}
and by \eqref{eq:cocycles_spin}
\begin{equation}\label{eq:action_U}
(\tau(f)x_{\tau^{-1}}^{1/2})(\omega)=f(\tau(\omega))\exp\Big[\frac12 \big(H(\omega)-H(\tau(\omega))\big)\Big].
\end{equation}
We also recall for each $\tau\in G$ an automorphism on $\mathfrak A$:
\[
u_\tau \pi_\mu(f)=U_\tau \pi_\mu(f)U_\tau^{*}.
\]
By the results of subsection \ref{subsec:conditional_expectation} we have a sequence of conditional expectations $(E_N)_N$:
\[
E_N(\pi_\mu(f))=\int_{G_N}u_\tau(\pi_\mu(f))d\tau.
\]
Now let us consider the projections $P_N=\int_{G_N}U_\tau d\tau$ and the vector $\Phi_N=P_N\Phi_\mu$. By \eqref{eq:U_again} and \eqref{eq:action_U}
\begin{eqnarray}\label{eq:dependent_rvs}
\Phi_N(\omega) &=&\int_{G_N}(U_\tau\Phi_\mu)(\omega) d\tau\nonumber\\
&=&\frac1{|G_N|}\sum_{\tau\in S_{\Lambda_N}}\exp\Big[\frac12 \big(H(\omega)-H(\tau(\omega))\big)\Big]. 
\end{eqnarray}
Defining unit vectors on $\mathcal H_\mu$ by 
\[
\Psi_N:=\frac1{\|\Phi_N\|}\Phi_N,
\]
we get vector states $\psi_N$ on $\mathfrak A$ by 
\begin{equation}\label{eq:vector_states}
\psi_N(f)=\langle \Psi_N, f\Psi_N\rangle_{\mathcal H_\mu},\quad f\in \mathfrak A.
\end{equation}
By Banach-Alaoglu theorem there is a subnet $(\psi_{N_k})_k$ such that it converges in weak$^*$-topology to a state, say $\psi$, on $\mathfrak A$. It is obvious that $\psi$ is $G$-invariant. We have the following result analogous to Theorem \ref{thm:martingale_convergence}.
\begin{thm}\label{thm:martingale_spin_system}
We assume that the state $\psi$ is faithful. Then the state $\psi$ is $\hat G$-invariant, i.e., $\psi(u_\tau(\pi_\mu(f)))=\psi(\pi_\mu(f))$ for all $f\in C(\Omega)$, and for any element $x\in \mathfrak A$, the sequence $(x_N)_N$, $x_N=E_N(x)$, is a decreasing martingale and has a limit $E_\infty(x):=\lim_N E_N(x)$.
\end{thm}
\Proof
The proof is exactly the same as in Theorem \ref{thm:martingale_convergence}.
\EndProof
\begin{rem}\label{rem:intuitive}
We were not able to show that the limit $\Phi_G:=\lim_N\Phi_N$ is non-zero. However, for the ferro magnetic Ising model, at least in the low-temparature regime, it must be true. The intuitive reasoning is as follows. Suppose the extreme case of zero temparature. Then the Gibbs measure has support on the configuration of all $+1$-spins or on the configuration of all $-1$-spins. In this case, obviously $x_\tau=1$ for all $\tau\in G$ (of course $\mu$ is $G$-invariant). Now if the temperature is not zero but sufficiently low, then there is a big cluster of same spins. In that case the exponent in the computation of $x_\tau$ in \eqref{eq:cocycles_spin} is near to zero for many $\tau$'s, meaning that $x_\tau$ is near $1$. So, the result of Propositon \ref{prop:sufficient_condition}. 
\end{rem}

\vskip 1 true cm
\noindent{\bf Acknowledgement}. We thank Mrs. Yoo Jin Cha for drawing the figure. A. Dhahri is a member of GNAMPA-INdAM and he has been supported by the MUR grant Dipartimento di Eccellenza 2023-2027 of Dipartimento di Matematica, Politecnico di Milano. The work of H. J. Yoo was supported by the National Research Foundation of Korea (NRF) grant funded by the Korean government (MSIT) (No. RS-2023-00244129).\\[1ex]

\noindent {\bf Conflict of Interest Statement}. There is no conflict of interest in this article.\\[1ex]
\noindent {\bf Data Availability Statement}. Data sharing not applicable to this article as no datasets were generated or analysed during the current study.

 \end{document}